\documentclass{svmult}

\usepackage{makeidx}         
\usepackage{graphicx}        
\usepackage{multicol}        
\usepackage[bottom]{footmisc}
\usepackage{url}
\usepackage{amsmath}
\usepackage{amssymb}
\usepackage{amsfonts}

\makeindex             

\begin{document}

\title*{A Personal Perspective on Numerical Analysis and Optimization}
 \titlerunning{History of Numerical Analysis}
\author{Desmond J. Higham\inst{1}}

 \institute{
School of Mathematics,
           University of Edinburgh,
           Edinburgh, EH9 3FD, UK}

%
%
\maketitle

\begin{abstract} {}
 I give a brief, non-technical, historical perspective on 
numerical analysis and optimization.
I also touch on emerging trends and 
future challenges.  
This content is based on the short presentation that I made at the  opening 
ceremony of 
\emph{The International Conference on Numerical Analysis and Optimization},
which was held at Sultan Qaboos University, Muscat, Oman, on 
January 6--9, 2020. 
Of course, the material covered here is necessarily incomplete and 
biased towards my own interests and comfort zones.
My aim is to 
give a feel for how the 
area has developed over the past few decades and how it may continue.
\end{abstract}

\section{Definitions}
\label{sec:defs}

Mathematicians love to make definitions. But
defining an area of mathematics is a thankless task.
The best one-sentence definitions that I can come up with 
for numerical analysis and optimization are as follows.
\begin{description}
\item[Numerical Analysis:] the design, analysis and implementation of
\emph{computational algorithms} to deliver approximate solutions to problems arising in \emph{applied mathematics}.
\item[Optimization:] the design, analysis and implementation of
\emph{computational algorithms} to approximate the 
\emph{best solution} 
to a problem arising in \emph{applied mathematics}
when there may be 
\emph{many feasible solutions}.
\end{description}
For alternative versions, I refer to the references
\cite{barrow-green15,Tref92,wright15}.

\section{Emergence}
\label{sec:em}
The unstoppable growth of interest in the use of computational algorithms  
can be attributed to two main factors.
First, technology has advanced rapidly and consistently since the 
digital computing era began in the 1950s.
The CDC6600, widely acknowledged to be the 
world's first ``supercomputer'', was introduced in 1964, achieving 
a  speed of 3 megaflops (that is, $3 \times {10}^{6}$ floating operations per second) \cite{H19}.
Today's fastest supercomputers can achieve 
petaflop speeds (${10}^{15}$ floating operations per second).
By contrast,
in his 1970 Turing Award Lecture \cite[1971]{wilk71a},
James Wilkinson discussed the use 
of mechanical desk calculators:
\begin{quote} \noindent
``It happened that some time after my arrival [at the National Physical
Laboratory in 1946],
a system of 18 equations arrived in Mathematics Division
and after talking around it for some time we finally decided to
abandon theorizing and to solve it \dots\
The operation was manned by Fox, Goodwin, Turing, and me,
and we decided on Gaussian elimination with complete pivoting.''
\end{quote}
Leslie Fox \cite{fox87}
noted that the computation referred to in this quotation
took about two weeks to complete.
By my estimation, this corresponds to around 
$0.003$ 
floating operations per second!

A second, and equally important, factor behind the rise of scientific computation 
is the availability of ever-increasing sources of data, caused by improvements in experimental techniques and, perhaps most notably, by the inexorable 
sensorization and digitization of our everyday lives.
Hence, although numerical analysis and optimization build on classical ideas 
that can be attributed to the likes of 
Newton, Euler, 
Lagrange and 
Gauss,
they continue to be shaped by current developments.    
Of course,
many disciplines 
make extensive use of 
computational techniques. For example,
 the word ``Computational'' often appears before the words 
Biology, Chemistry, Physics, and Social Science.
Furthermore, Computational Science and Engineering 
\cite{Retal18}
is a well established discipline
that is often referred to as the 
third leg
of the 
science and engineering stool, equal in stature to 
 observation and theory.
 In addition,
 many graduate schools  
 now offer courses with titles such as 
  ``Data Analytics''  
  and ``Data Science.''
  Although there is clearly much overlap, my view is that 
  numerical analysis and optimization 
   have a distinct
  role of focusing on the \emph{design} and \emph{analysis} of
   algorithms for problems in applied mathematics, 
   in terms of complexity, accuracy and stability,   
   and hence they are
    informed by, but not driven by, application fields.

\section{Reflections}
\label{sec:refs}

My own exposure to numerical analysis and optimization dates back to the mid 1980s when I enrolled 
on an MSc course on ``Numerical Analysis and Programming'' at the University of Manchester. The course, in which optimization was treated as a branch of numerical analysis, made heavy use of the series of textbooks
published by Wiley that was written by members of the 
highly influential numerical analysis group at the University of Dundee  
\cite{GAW09}.
Here are 
the topics, and most recent versions of these books: approximation theory \cite{W80}, 
numerical methods for ODEs \cite{Lam91},
numerical methods for PDEs \cite{MG80,WM85} and 
optimization \cite{F2000}.
An 
additional topic was numerical linear algebra
\cite{GVLbook}.
Each of these topics remains active and highly relevant.
Numerical linear algebra and optimization are often important 
building blocks 
within larger computational tasks,
and hence their popularity has never waned.
Partial differential equations lie at the heart of most models in the natural and engineering sciences, and they come in many varieties, giving rise to 
an ever-expanding problem set.
Timestepping methods for ODEs 
gained impetus through the concept of geometric integration 
\cite{HLW06}
and now 
play a prominent role
in the  
development of 
tools for statistical sampling \cite{mcmc11}. 
ODE simulation also forms a key component in certain classes of neural network, 
as described in 
\cite{NIPS2018_7892},
which received a 
Best Paper Award at NeurIPS 2018, a leading conference in machine learning.
Approximation theory
lies at the heart of 
the current deep learning revolution 
\cite{HHdl2019,CBG15},
and, in particular,  
understanding very high dimensional data
spaces and/or parameter spaces remains a fundamental challenge.

\section{Impact}
\label{sec:imp}

In a special issue of the journal 
Computing in Science and Engineering,
Jack Dongarra and Francis Sullivan
published their top-ten list of 
algorithms that had the 
``greatest influence on the development and practice of science and engineering in the 20th century''
\cite{DS2000}. These were
\begin{itemize}
\item
Metropolis Algorithm for Monte Carlo
\item 
Simplex Method for Linear Programming
\item 
Krylov Subspace Iteration Methods
\item 
The Decompositional Approach to Matrix Computations
\item 
The Fortran Optimizing Compiler
\item 
QR Algorithm for Computing Eigenvalues
\item 
Quicksort Algorithm for Sorting
\item 
Fast Fourier Transform
\item 
Integer Relation Detection
\item 
Fast Multipole Method
\end{itemize}
Here, the word ``algorithm'' is being used in a very general sense, but it is clear 
that 
most of these
achievements have 
ideas from numerical analysis and optimization at their heart.

Researchers in numerical analysis and optimization
are in the advantageous position that their work is not only recorded in
journals and textbooks, but may also be made 
operational through public domain software.
Many authors now deposit code alongside their academic publications,
and
state-of-the-art code is available in 
 a wide range of  
languages and platforms, including 
FORTRAN, C, R, MATLAB,  Maple,  
(Scientific) Python and the more recent Julia \cite{Julia17}.

\section{Momentum}
\label{sec:mom}

Judging by the level of activity around graduate classes, 
seminar series, conferences and journals, 
there is a strong pull for further research 
in numerical analysis and optimization.
Particularly active, and overlapping, directions include
\begin{itemize}
\item dealing with, or exploiting, randomness, both in the simulation of  
mathematical models that are inherently stochastic   
\cite{LPS14}
and in the use of 
randomization in solving deterministic problems \cite{CHM20,MT20},
\item  accurately and efficiently simulating mathematical models that operate over a vast range of temporal or spatial scales
\cite{E2011,ecm2017},
\item tackling problems of extremely high dimension, notably
 large-scale optimization and inverse problems in machine learning
 and imaging 
\cite{AMOS19,BCN18},
\item 
exploiting the latest computer architectures and 
designing algorithms that efficiently trade off between 
issues 
such as 
 memory bandwidth, data access, communication 
and, perhaps most topically, the use of low precision special function units
\cite{aabc20}.
\end{itemize}

Such items, and may others, emphasize that important challenges remain 
for 
researchers in numerical analysis and optimization  
in the design, evaluation and extension of 
the modern computational scientist's toolbox.

\begin{acknowledgement}
 The author is 
 supported by 
  Programme Grant EP/P020720/1
  from the Engineering and Physical Sciences Research Council of the UK.
\end{acknowledgement}

\bibliographystyle{siam}
\bibliography{hawkesrefs,pcam}

\end{document}